\documentclass[11pt]{amsart}

\usepackage{latexsym, graphics, psfrag, amscd, amssymb}
\newtheorem{thm}{Theorem}
\newtheorem{tm}[thm]{Theorem}
\newtheorem{fact}[thm]{Fact}
\newtheorem{lm}[thm]{Lemma}
\newtheorem{rem}[thm]{Remark}
\newtheorem{cor}[thm]{Corollary}
\newtheorem{pro}[thm]{Proposition}
\newtheorem{ex}[thm]{Example}
\newtheorem{df}[thm]{Definition}

\newcommand{\Si}{\Sigma}

\newcommand{\cM}{\mathcal{M}}

\newcommand{\cC}{{\mathcal{C}}}
\newcommand{\cD}{{\mathcal{D}}}
\newcommand{\Z}{\mathbb{Z}}
\newcommand{\R}{\mathbb{R}}
\newcommand{\C}{\mathbb{C}}
\newcommand{\RP}{\mathbb{RP}}

\newcommand{\Hol}{\mathrm{Hol}}
\newcommand{\Ker}{\mathrm{Ker}}

\newcommand{\Hom}{\mathrm{Hom}}

\newcommand{\Prin}{\mathrm{Prin}}
\newcommand{\im}{\mathrm{Im}}

\newcommand{\fg}{\mathfrak{g}}
\newcommand{\ft}{\mathfrak{t}}
\newcommand{\fh}{\mathfrak{h}}

\newcommand{\tG}{ \tilde{G} }

\newcommand{\tg}{\tilde{g}}
\newcommand{\ta}{\tilde{a}}
\newcommand{\tb}{\tilde{b}}
\newcommand{\tc}{\tilde{c}}

\newcommand{\td}{\tilde{d}}
\newcommand{\tT}{\tilde{T}}

\newcommand{\ra}{\rightarrow}
\newcommand{\hol}{\Hom(\pi_1(\Si),G)}
\newcommand{\holC}{\Hom_{\Cr}(\pi_1(\Si),G)}
\newcommand{\ab}{a_1,b_1,\ldots, a_\ell,b_\ell}

\newcommand{\slr}{\Sigma^{\ell,r}_0}
\newcommand{\sone}{\Sigma^{\ell,r}_1}
\newcommand{\stwo}{\Sigma^{\ell,r}_2}
\newcommand{\lri}{\Sigma^{\ell,r}_i}
\newcommand{\Cr}{\cC_1,\ldots,\cC_r}
\newcommand{\Ct}{\cC_1\times\cdots\times\cC_r}

\newcommand{\mlri}{\cM_G(\lri;\Cr)}

\newcommand{\Xz}{X_G^{\ell,r,0}(\Cr;z)}
\newcommand{\Xo}{X_G^{\ell,r,1}(\Cr;z)}
\newcommand{\Xt}{X_G^{\ell,r,2}(\Cr;z)}

\begin{document}

\parskip=0.35\baselineskip
\baselineskip=1.2\baselineskip

\title{Connected Components of Spaces of Surface Group Representations II}

\author{Nan-Kuo Ho}
\address{The Fields Institute, Toronto, Canada }
\address{ Department of
Mathematics, National Cheng-Kung University, Taiwan}
\email{nho@fields.utoronto.ca \\ nankuo@mail.ncku.edu.tw}

\author{Chiu-Chu Melissa Liu}
\address{Department of Mathematics \\ Harvard University}
\email{ccliu@math.harvard.edu}


\subjclass{53}
\date{September 21, 2004}

\begin{abstract}
In \cite{HL1}, we discussed the connected
components of the space of surface group representations for any
compact connected semisimple Lie group and any closed
compact (orientable or nonorientable) surface. In this sequel,
we generalize the results in \cite{HL1} in two directions:
we consider general compact connected Lie groups, and we
consider all compact surfaces, including the ones
with boundaries. We also interpret our results in terms
of moduli spaces of flat connections over compact surfaces.
\end{abstract}
\maketitle

\section{introduction}
In \cite{G2}, W. M. Goldman computed the number of connected components of
spaces of surface group representations $\hol$ for $G=PSL(2,\C)$
or $PSL(2,\R)$, and $\Si$ is a Riemann surface of genus bigger
than $1$. In the same paper, he also made a conjecture for the
general connected complex semisimple Lie group $G$ that there is a
bijection between $\pi_0(\hol/G)$ and
$H^2(\Si,\pi_1(G))\cong\pi_1(G)$, where the $G$-action on $\hol$ is induced
by the conjugate action of $G$ on itself. This conjecture was proved by J. Li in \cite{li}.

In \cite{HL1}, we computed $\pi_0(\hol/G)$ for any compact connected
semisimple Lie group $G$ and any closed (orientable or nonorientable) surface
$\Si$ with negative Euler characteristic, except for the connected sum of a
torus and a Klein bottle.  The first goal of this note is to generalize
the results in \cite{HL1} to general compact connected Lie groups.
The result for orientable surfaces is known:

\begin{tm}\label{thm:irzero}
Let $\Si$ be a connected, closed, compact, orientable surface with genus $\ell>0$.
Let $G$ be a compact connected Lie group, and let $G_{ss}=[G,G]$ be the
maximal connected semi-simple subgroup of $G$.
Then there is a bijection
$$
\pi_0(\hol/G)\to \pi_1(G_{ss}).
$$
\end{tm}
The group $\pi_1(G_{ss})$ in Theorem \ref{thm:irzero} is a finite abelian group.
As we will explain in Section \ref{sec:moduli}, a proof of Theorem
\ref{thm:irzero} can be extracted from \cite{ym}. In this note,
we will give another proof of Theorem \ref{thm:irzero}, and
derive the following result for non-orientable surfaces:

\begin{tm}\label{thm:rzero}
Let $\Si$ be a closed, compact, nonorientable surface which is
homeomorphic to the connected sum of $m$ copies of the real
projective plane, where $m\neq 1,2,4$. Let $G$ be a compact
connected Lie group. Then there is a bijection
$$
\pi_0(\hol/G)\to \pi_1(G)/2\pi_1(G)
$$
where $2\pi_1(G)$ denote the subgroup
$\{k^2\mid k\in\pi_1(G)\}$ of the abelian group $\pi_1(G)$.
\end{tm}

The group $\pi_1(G)/2\pi_1(G)$ in Theorem \ref{thm:rzero} fits in
the following short exact sequence of abelian groups:
$$
1\to \pi_1(G_{ss})/2\pi_1(G_{ss})\to \pi_1(G)/2\pi_2(G)\to
(\Z/2\Z)^{\dim H}\to 1
$$
where $G_{ss}=[G,G]$ and $H$ is the connected component of the identity
of the center of $G$. In particular, $\pi_1(G)/2\pi_1(G)$ is a finite
abelian group.

Theorem \ref{thm:irzero} and Theorem \ref{thm:rzero}
generalize the case where $G$ is semisimple considered in
\cite{HL1} and the case where $G_{ss}$ is simply connected
considered in \cite{HL2}.

The second goal of this note is to generalize the above results
to compact surfaces with boundaries. The same question becomes trivial
in this case:  $\pi_1(\Si)$ is a free group when $\Si$ is a compact
surface with nonempty boundary, so $\hol$ and $\hol/G$ are connected
for any compact connected Lie group $G$. However, the question
becomes interesting when we impose some boundary conditions.
Let $\Si$ be a (possibly nonorientable) compact surface with boundary
components $B_1,\ldots, B_r$. Let
\begin{equation}\label{eqn:holC}
\holC=
\{\phi\in \hol\mid \phi([B_j])\in \cC_j\}
\end{equation}
where $[B_j]\in \pi_1(\Si)$ and $\Cr$ are conjugacy classes in $G$, also known as
{\em markings} in this context. Note that (\ref{eqn:holC})
reduces to $\hol$ when $r=0$. The conjugate action of $G$ on itself
induces a $G$-action on $\holC$. We will compute
$$
\pi_0(\holC/G)=\pi_0(\holC).
$$

For orientable surfaces, we have
\begin{tm}\label{thm:izero}
Let $\Si$ be a connected, closed, compact, orientable surface with $\ell>0$
handles and $r>0$ boundary components.
Let $G$ be a compact connected Lie group, and let $G_{ss}=[G,G]$ be the
maximal connected semi-simple subgroup of $G$.
Then
$$
\holC/G
$$
is nonempty iff $d_1\cdots d_r\in G_{ss}$ for some
(and therefore for all) $(d_1,\ldots,d_r)\in \Ct$.
For any conjugacy classes $\Cr$ of $G$ such that $d_1\cdots d_r\in G_{ss}$
for some $(d_1,\ldots, d_r)\in\Ct$, there is a bijection
$$
\pi_0(\holC/G)\to \pi_1(G_{ss})/J_{\Cr}
$$
where $J_{\Cr}$ is a subgroup of $\pi_1(G_{ss})$ defined
in Section \ref{sec:conjugate}.
\end{tm}

Let $\Si$ and $G$ be as Theorem \ref{thm:izero}. 
By Theorem \ref{thm:izero}, when $G$ is semisimple and simply connnected,
for example, $SU(n)\,(n\geq 2)$, $Spin(n)\,(n\geq 3)$, $Sp(n)\,(n\geq 1)$,
$\holC/G$ is nonempty and connected  for any conjugacy classes $\Cr$; 
when $G$ is not semisimple and $G_{ss}$ 
is simply connected, for example, $U(n)$ or compact torus, $\holC/G$ is either 
empty or connected. When $G=SO(n)\, (n\geq 3)$, $\holC/G$ is nonempty 
and has either one or two connected components. We will see later that for 
generic conjugacy classes $\Cr$ in $SO(n)$, $J_{\Cr}$ is trivial, and in 
that case $\holC/G$ has two connected components.
 
For nonorientable surfaces, we have
\begin{tm}\label{thm:nonzero}
Let $\Si$ be a connected, compact, nonorientable surface which is
homeomorphic to the connected sum of a sphere with $r$ holes and
$m$ copies of the real projective plane, where $m\neq 1,2,4$.
Let $G$ be a compact connected Lie group. Then there is a bijection
$$
\pi_0(\holC/G)\to \pi_1(G)/J'_{\Cr}
$$
where $J'_{\Cr}$ is the subgroup of $\pi_1(G)$ generated by
$2\pi_1(G)$ and $J_{\Cr}$.
\end{tm}
 
Note that $\pi_1(G)/J'_{\Cr}$ is the quotient group of the finite
abelian group $\pi_1(G)/2\pi_1(G)$. Let $\Si$ and $G$ be
as in Theorem \ref{thm:nonzero}. When $G$ is semisimple and simply 
connected, for example,
 $SU(n)\, (n\geq 2)$, $Spin(n)\,(n\geq 3)$, $Sp(n)\,(n\geq 1)$,
$\hol/G$ is nonempty and connected for any conjugacy classes
$\Cr$ in $G$. When $G=U(n)$, $\pi_1(G)=\Z$ and 
$J_{\Cr}\subset \pi_1(G_{ss})$ is
trivial, so $\holC/G$ has two connected components.
When  $G=SO(n)$, we have  
$$
\pi_1(G)/J'_{\Cr}=\pi_1(G_{ss})/J_{\Cr}
$$
since $\pi_1(G_{ss})=\pi_1(G)=\Z/2\Z$ and $2\pi_1(G)$ is trivial, so 
the answer is the same as in the orientable case.

Our proofs of the above results rely on a result of Alekseev, Malkin, and Meinrenken
on group-valued moment maps (Fact \ref{thm:fiber}).

The rest of this paper is organized as follows. In Section \ref{sec:pre}, we give some
preliminaries on the structure of compact Lie groups and introduce some
definitions. In Section \ref{sec:special}, we derive
Theorem \ref{thm:irzero}--\ref{thm:nonzero} under the additional assumption that $G$
is simply connected. The general case is treated in Section
\ref{sec:general}: the results for orientable surfaces (Theorem
\ref{thm:irzero} and Theorem \ref{thm:izero}) are proved in Section
\ref{sec:orientable}, and the results for nonorientable surfaces (Theorem
\ref{thm:rzero} and Theorem \ref{thm:nonzero}) are proved in Section \ref{sec:nonorientable}.
In Section \ref{sec:moduli}, we give a geometric interpretation of
the above results to the moduli spaces of flat bundles over compact surfaces.

\section{Preliminaries}\label{sec:pre}

\subsection{Compact connected Lie groups}\label{sec:compact}
Let $G$ be a compact connected Lie group. Let $G_{ss}=[G,G]$ be its
commutator group. Then $G_{ss}$ is the maximal connected semisimple
subgroup of $G$. Let $H$ be the connected component of the identity
of the center $Z(G)$ of $G$. Then $H$ is a compact torus.
The map $\phi: H\times G_{ss}\to G=HG_{ss}$
given by $(h,g)\mapsto hg$ is a finite cover which is also a group
homomorphism. The kernel of $\phi$ is isomorphic to $D=H\cap G_{ss}\subset Z(G_{ss})$,
which is a finite abelian group. Note that $G_{ss}$ is a normal subgroup
of $G$, and the quotient
$$
G/G_{ss}\cong H/D
$$
is a compact torus.

Let $\rho_{ss}:\tG_{ss}\to G_{ss}$ be the
universal covering map which is also a group homomorphism. Then $\tG_{ss}$ is
a compact, connected, simply connected Lie group, and
$\Ker(\rho_{ss})$ is a subgroup of $Z(\tG_{ss})$.

Let $\fg$ and $\fh$ be the Lie algebras of $G$ and $H$ respectively,
and let $\exp_H:\fh\to H$ be the exponential map. The map
$$
\rho: \tG=\fh\times \tG_{ss} \to G
$$
given by
$$
(X,g)\to \exp_H(X)\rho_{ss}(g)
$$
is the universal covering map which is also a group homomorphism.
Here we give a group structure for $\fh$ the addition operation in
the vector space $\fh$ and $0$ as the identity element for the
group structure. Notice that $\tG$ is not compact.
 Let $\pi_1: H\to H/D$ and $\pi_2: G_{ss}\to G_{ss}/D$ be
natural projections. Then
$$
\pi_1\circ \exp_H :\fh\to H/D\cong G/G_{ss}
$$
is the universal covering map which is also a group homomorphism, and
$$
\check{\Lambda}=\Ker(\pi_1\circ\exp_H)\cong \pi_1(G/G_{ss})\cong \Z^{\dim H}.
$$
We have
\begin{eqnarray*}
\Ker(\rho) &=& \{ (X,g)\in \fh\times \tG_{ss}\mid \exp_H(X)\rho_{ss}(g)=e\}\\
&& \subset \check{\Lambda} \times \Ker(\pi_2\circ\rho_{ss})\subset
\check{\Lambda}\times Z(\tG_{ss}) \subset \fh\times Z(\tG_{ss})=Z(\tG).
\end{eqnarray*}

The map $(X,g)\mapsto X$ defines a surjective group homomorphism
$p:\Ker(\rho)\to \check{\Lambda}$. The kernel of $p$ is
$\{0\}\times \Ker(\rho_{ss})\cong \Ker(\rho_{ss})$. So
we have an exact  sequence of abelian groups
\begin{equation}\label{eqn:kernel}
1\to \Ker(\rho_{ss})\to \Ker(\rho)\to \check{\Lambda}\to 1,
\end{equation}
which can be rewritten as
\begin{equation}
1\to \pi_1(G_{ss})\to \pi_1(G)\to \pi_1(G/G_{ss})\to 1.
\end{equation}
where $\pi_1(G_{ss})$ is a finite abelian group, and
$$
\pi_1(G/G_{ss})=\pi_1(H/D)\cong \Z^{\dim H}.
$$

\subsection{Conjugacy classes}\label{sec:conjugate}
Let $G$ be a connected Lie group with center $Z(G)$.
Given $a\in G$, let $G\cdot a$ denote the orbit of $a$ of the
conjugate action of $G$ on itself. The left multiplication
by $z\in Z(G)$ gives a bijection $G\cdot a \to G\cdot (za)$.
Given a conjugacy class $\cC=G\cdot a$ and $z\in Z(G)$, let
$z\cC$ denote the conjugacy class $G\cdot (za)$.

Let $K$ be a subgroup of $Z(G)$. Let $K$ acts on $G$ by left multiplication.
This $K$-action on $G$ commutes with the conjugate action of $G$ and induces
a $K$-action on $Con(G)$, the set of conjugacy classes of $G$.
Given a conjugacy $\cC\in Con(G)$, let $K_{\cC}$ denote the stabilizer
of $\cC$ under the $K$-action on $Con(G)$. It is straightforward to check that
$K_{\cC}=K_{z\cC}$ for any $\cC\in Con(G)$ and $z\in Z(G)$.

Now let $G$ be a compact connected Lie group, and let
$\rho:\tG\to G$ be as in Section \ref{sec:compact}.
Then $\tG$ is a noncompact connected Lie group. 
The following is obviously true.
\begin{lm}\label{thm:ss}
Let $K=\Ker\rho\subset Z(\tG)=\fh\times Z(\tG_{ss})$, and let
$\tilde{K}=\Ker(\rho_{ss})\subset Z(\tG_{ss})$. Then any 
conjucay class $\cD\in Con(\tG)$ is of the form
$$
\cD=\{X\}\times \tilde{\cC}
$$
for some $X\in \fh$ and $\tilde{\cC}\in Con(\tG_{ss})$, and we have
$$ 
K_\cD= \{0\}\times \tilde{K}_{\tilde{\cC}}
\subset \{0\}\times \Ker(\rho_{ss}).
$$
\end{lm}

Given a conjugacy class $\cC$ in $G$, each connected
component of $\rho^{-1}(\cC)$ is a conjugacy class of $\tG$. Let
$\cD$ be a connected component of $\rho^{-1}(\cC)$. Then
$\cD\to \cC$ is a covering map of finite degree. Actually,
the degree of $\cD\to \cC$ is equal to the cardinality of
the subgroup $K_\cD$ of $\{0\}\times \Ker(\rho_{ss})\cong \pi_1(G_{ss})$.
Let
$$
J_{\Cr}\subset \{0\}\times \Ker(\rho_{ss})\cong \pi_1(G_{ss})
$$
be the subgroup generated by $K_{\cD_1}, \ldots, K_{\cD_r}$, where
$\cD_j$ is a connected component of $\rho^{-1}(\cC_j)$. Note that
the definition of $J_{\Cr}$ does not depend on the choices of $\cD_1,\ldots,\cD_r$
because any connected component of $\rho^{-1}(\cC_j)$ is of the
form $z\cD_j$ for some $z\in Z(\tG)$ and $K_{z\cD_j}= K_{\cD_j}$.

\begin{ex}\label{ex:so}
$G=SO(3)=G_{ss}$, $\tG=Spin(3)=\tG_{ss}$ ,and $\rho_{ss}=\rho:\tG\ra G$. 
We view $Spin(3)$ as a subset of the real Clifford algebra $C_3$. 
Choose a maximal torus $T$ of $SO(3)$ as follows:  
$$
T=\{
R_\theta=\left(\begin{array}{ccc}
 \cos\theta & -\sin\theta & 0\\
\sin\theta & \cos\theta & 0\\
 0 & 0 & 1
\end{array}\right)\mid  -\pi\leq \theta \leq \pi \}\cong U(1)
$$
Then
$$
\tT=\rho^{-1}(T)=\{\eta_\phi=\cos\phi-\sin\phi e_1 e_2\mid -\pi\leq
\phi\leq \pi\} \cong U(1)
$$
is a maximal torus of $Spin(3)$. We have $\rho(\eta_\phi)=R_{2\phi}$, 
so $\rho|_{\tT}:\tT\to T$ can be identified with 
$$
s:U(1)\to U(1),\ \ u\mapsto u^2.
$$

Let $W=\{\pm 1\}$ be the Weyl group of $SO(3)$ which is also the Weyl
group of $Spin(3)$. Then $W$ acts on $SO(3)$ and $Spin(3)$ by
$(-1)\cdot R_\theta = R_{-\theta}$ and $(-1)\cdot \eta_\phi=\eta_{-\phi}$,
respectively. Let $\cC_\theta \in Con(SO(3))$ be the conjugacy
class represented by $R_\theta \in T\subset SO(3)$, and let
$\cD_\phi \in Con(Spin(3))$ be the conjugacy class represented
by $\eta_\phi \in\tT\subset Spin(3)$. Then
\begin{eqnarray*}
Con(SO(3))=T/W
&=&\{\cC_\theta \mid 0\leq \theta \leq \pi\}\cong [0,\pi]\\
Con(Spin(3))=\tT/W
&=&\{\cD_\phi\mid 0\leq \theta \leq \pi \}\cong [0,\pi]
\end{eqnarray*}
$\rho|_{\tT}:\tT\to T$ induces a map $\hat{\rho}: \tT/W \to T/W$ which
can be identified with
$$
\hat{s}:[0,\pi]\to [0,\pi],\ \ 
\phi \mapsto\left\{ \begin{array}{ll} 
2\phi & 0\leq \phi \leq \frac{\pi}{2}\\
2(\pi-\phi)& \frac{\pi}{2} \leq \phi \leq \pi.
\end{array}\right.
$$
So for $\theta \in [0,\pi]$, we have
$\hat{\rho}^{-1}(\cC_\theta)=\{\cD_{\frac{\theta}{2}},
\cD_{\pi-\frac{\theta}{2}}\}$. Note that $\hat{\rho}^{-1}(\cC_\theta)$
consists of two conjugacy classes unless $\cC_\theta=\cC_\pi$, which
is represented by 
$$
R_\pi=\left(\begin{array}{ccc}
-1 & 0 & 0\\
0 & -1 & 0\\
0 & 0 & 1
\end{array}
\right).
$$
In this case, $\hat{\rho}^{-1}(\cC_\pi)$ consists of a single
conjugacy class $\cD_{\frac{\pi}{2}}$.
So for $\theta \in [0,\pi], \theta\neq \pi$,
$\rho^{-1}(\cC_\theta)=\cD_{\frac{\theta}{2}}\cup
\cD_{\pi-\frac{\theta}{2}} \to \cC_\theta$ is the trivial
(disconnected) double cover, while
$\rho^{-1}(\cC_{\pi})=\cD_{\frac{\pi}{2}}\to \cC_\pi$ is 
a connected nontrivial double cover.

Let $K=\{\pm 1\}=\Ker(\rho_{ss})=\Ker(\rho)$. Then $K$ acts on 
$\tT$ by $(-1)\cdot \eta_\phi = -\eta_\phi = \eta_{\phi-\pi}$.
So $K$ acts on $Con(Spin(3))$ by $(-1)\cdot \cD_\phi = 
\cD_{\phi-\pi} = \cD_{\pi-\phi}$, i.e., 
$(-1)\cdot \phi =\pi-\phi$ under the identification
$Con(Spin(n))\cong [0,\pi]$. It is clear from the
the above explicit description that for $\phi\in[0,\pi]$, 
$$
K_{D_\phi}=\left\{\begin{array}{ll}
\{    1\},& \textup{ if  }\phi \neq\frac{\pi}{2},\\
\{\pm 1\},& \textup{ if  }\phi =\frac{\pi}{2}.
\end{array}\right.
$$
Then for
$\Cr\in Con(SO(3))$,
$$
J_{\Cr}
=\left\{\begin{array}{ll}
\{    1\},& \textup{ if none of }\Cr\textup{ is equal to }\cC_{\pi},\\
\{\pm 1\},& \textup{ if one of  }\Cr\textup{ is equal to }\cC_{\pi}.
\end{array}\right.
$$

\end{ex}

The analysis in Example \ref{ex:so} can be generalized to show that
$J_{\Cr}$ is trivial for generic conjugacy classes 
$\Cr\in Con(SO(n))$, $n\geq 3$.

\subsection{Representation varieties}
\label{sec:variety}

\begin{df}\label{def:variety}
Let $G$ be a connected Lie group.
Given nonnegative integers $\ell, r$ such that $(\ell,r)\neq(0,0)$ and conjugacy
classes $\Cr\in Con(G)$, define
\begin{eqnarray*}
\Phi_{G,\ell,r}^{\Cr}: && G^{2\ell}\times\Ct \longrightarrow G\\
&&(a_1,b_1,\ldots,a_\ell,b_\ell,d_1,\ldots,d_r)
\mapsto \prod_{i=1}^\ell [a_i,b_i]\prod_{j=1}^r d_j.
\end{eqnarray*}
Let $\Phi_{G,0,0}:\{e\}\to G$ be the inclusion of the trivial subgroup.

Given nonnegative integers $\ell,r$, conjugacy classes $\Cr\in Con(G)$,
and $z\in Z(G)$, define
\begin{eqnarray*}
&&\Xz= \left(\Phi_{G,\ell,r}^{\Cr}\right)^{-1}(z)\\
& & \Xo\\
&=&\{ (a_1,b_1,\ldots,a_\ell,b_\ell,d_1,\ldots,d_r,c)
      \in G^{2\ell}\times\Ct\times G \mid
      \prod_{i=1}^\ell[a_i,b_i] \prod_{j=1}^r d_j c^2 =z\}\\
& & \Xt\\
&=&\{ (a_1,b_1,\ldots,a_\ell,b_\ell,d_1,\ldots,d_r,c_1,c_2)
      \in G^{2\ell}\times\Ct\times G^2 \mid
      \prod_{i=1}^\ell[a_i,b_i] \prod_{j=1}^r d_j c_1^2 c_2^2 =z\}
\end{eqnarray*}
\end{df}

Since $G$ is a Lie group and any conjugacy class of $G$ is a homogeneous
space, for any nonnegative integers $\ell,r$ and $i=1,2,3$
$$
G^{2\ell}\times \Ct \times G^i
$$
is a manifold (it is a point when $\ell=r=i=0$). In this topology,
$X^{\ell,r,i}_G(\Cr;z)$ is a closed subset.
Let $X^{\ell,r,i}_G(\Cr;z)$ be equipped with the induced topology
as a closed subset of $X^{\ell,r,i}_G(\Cr;z)$.

Let $\slr$ be the compact, connected, orientable surface with $\ell$ handles and $r$
boundary components $B_1,\ldots, B_r$, where $\ell, r$ are nonnegative integers.
Let $\sone$ be the connected sum of $\slr$ and $\RP^2$, and let
$\stwo$ be the connected sum of $\slr$ and a Klein bottle.
Any compact surface is of the form $\lri$, where $\ell,r$ are nonnegative
integers and $i=0,1,2$. It is orientable iff $i=0$; it is closed iff $r=0$.
Then $\Hom_{\Cr}(\pi_1(\lri),G)$ can be identified with $X^{\ell,r,i}_G(\Cr;e)$.

\section{Simply connected case} \label{sec:special}

In this section, we consider a compact, connected, simply connected Lie group
$G$. In particular, $G$ is semisimple.

\begin{fact}[{\cite{dhm}}]\label{thm:fiber}
Let $G$ be a compact, connected, simply connected Lie group.
Let $\ell$ be a positive integer, and let $\Cr$ be any conjugacy
classes in $G$. Then for any $g\in G$, $(\Phi_{G,\ell,r}^{\Cr})^{-1}(g)$
is nonempty and connected.
\end{fact}

\begin{rem}
The reason of Fact \ref{thm:fiber} is that $\{ G^{2\ell}\times\Ct, \Phi_{G,\ell,r}^{\Cr}\}$
is a quasi-Hamiltonian system and thus the preimage of the
moment map $\Phi_{G,\ell,r}^{\Cr}$ is nonempty and connected under the above assumption of $G$.
\end{rem}

We will derive the following result from Fact \ref{thm:fiber}.
\begin{pro}\label{thm:connected}
Let $G$ be a compact, connected, simply connected Lie group.
Let $\ell$ be a positive integer, let $r$ be a nonnegative integer,
let $\cC_1,\ldots,\cC_r$ be any conjugacy classes in $G$, and let
$z\in Z(G)$.  Then
$$
\Xz,\ \ \Xo
$$
are nonempty and connected;
$$
\Xt
$$
is nonempty, and is connected if $\ell>1$.
\end{pro}

\paragraph{\bf Proof}
By Fact \ref{thm:fiber},
$$
\Xz=\left(\Phi_{G,\ell,r}^{\Cr}\right)^{-1}(z).
$$
is nonempty and connected. In the rest of the proof,
we will consider $\Xo$ and $\Xt$.

We fix a maximal torus $T$ in $G$.
For $r=0$, let $s\in T$ be a square root of $z$. For $r>0$,
$T\cap \cC_j$ is nonempty for
$j=1,\ldots,r$. We fix $t_j\in T\cap \cC_j$, and let
$s\in T$ be a square root of $(t_1\cdots t_r)^{-1}z$.

Consider maps
\begin{eqnarray*}
Q_1:X_G^{\ell,r,1}(\Cr;z) &\to& G\\
(a_1,b_1,\ldots,a_\ell,b_\ell,d_1,\ldots,d_r, c)&\mapsto& c\\
Q_2:X_G^{\ell,r,2}(\Cr;z) &\to& G^2\\
(a_1,b_1,\ldots,a_\ell,b_\ell,d_1,\ldots,d_r, c_1,c_2)&\mapsto& (c_1,c_2)
\end{eqnarray*}
By Fact \ref{thm:fiber}, $Q_1^{-1}(c)$ is nonempty and connected for any $c\in
G$, and $Q_2^{-1}(c_1,c_2)$ is nonempty and connected for any $(c_1,c_2)\in G^2$.
So  $X_G^{\ell,r,1}(\Cr;z)$ and $X_G^{\ell,r,2}(\Cr;z)$ are  nonempty. We
will show that
\begin{enumerate}
\item For any $c_0\in G$ there is a path
      $\gamma_1:[0,1]\to X_G^{\ell,r,1}(\Cr;z)$  such that
      $\gamma_1(0)\in Q^{-1}_1(s)$ and $\gamma_1(1)\in Q^{-1}_1(c_0)$.
\item If $\ell>1$, then for any $(c_1,c_2)\in G^2$ there is a path
      $\gamma_2:[0,1]\to X_G^{\ell,r,2}(\Cr;z)$  such that
      $\gamma_2(0)\in Q_2^{-1}(s,e)$ and $\gamma_2(1)\in Q_2^{-1}(c_1,c_2)$.
\end{enumerate}
This will complete the proof.

Given $c_0, c_1, c_2\in G$,  there exist $g_0,g_1,g_2\in G$ such that
$g_i^{-1} c_i g_i\in T$. Let $\fg$ and $\ft$ be the Lie algebras of $G$ and $T$,
respectively. Let $\exp:\fg\to G$ be the exponential map.
Then there exist $\xi_0,\xi_1,\xi_2\in\ft$ such that
$$
s^{-1}g_0^{-1} c_0 g_0=\exp (\xi_0),\ \
s^{-1}g_1^{-1} c_1 g_1=\exp (\xi_1),\ \
g_2^{-1} c_2 g_2=\exp (\xi_2).
$$
Let $W$ be the Weyl group of $G$, and let $w\in W$ be a Coxeter element (cf:\cite{Hu}).
The linear map $w:\ft\to \ft$ has no eigenvalue equal to $1$,
so there exists $\xi_i'\in\ft$ such that $w\cdot \xi_i'-\xi_i'=\xi_i$.
Recall that $W=N(T)/T$, where $N(T)$ is the normalizer of $T$ in $G$,
so $w=aT\in N(T)/T$ for some $a\in G$. We have
$$
a\exp(t\xi_i')a^{-1} \exp(-t\xi_i')=\exp(t\xi_i)
$$
for any $t\in \R$.

The group $G$ is connected, so there exists a path $\tg_i:[0,1]\to G$ such that
$\tg_i(0)=e$ and $\tg_i(1)=g_i$, where $i=0,1,2$. Define $\gamma_1:[0,1]\longrightarrow
G^{2\ell}\times \Ct \times G$ by
$$
\gamma_1(t)=(a(t), b(t), e,\ldots,e, d_1(t),\ldots,d_r(t), c(t)),
$$
where
\begin{eqnarray*}
a(t)&=&\tg_0(t) a \tg_0(t)^{-1}, \\
b(t)&=&\tg_0(t)\exp(-2t\xi_0')\tg_0(t)^{-1}, \\
d_j(t)&=& \tg_0(t) t_j \tg_0(t)^{-1},\ \ j=1,\ldots,r\\
c(t)&=& \tg_0(t) s \exp(t\xi_0)\tg_0(t)^{-1}.
\end{eqnarray*}
Then the image of $\gamma_1$ lies in $X_G^{\ell,r,1}(\Cr;z)$,
$$
\gamma_1(0)=(a,e,e,\ldots,e,t_1,\ldots,t_r,s)\in Q_1^{-1}(s),
$$
and
\begin{eqnarray*}
\gamma_1(1)&=&(g_0ag_0^{-1}, g_0\exp(-2\xi_0')g_0^{-1}, e,\ldots,e,
 g_0 t_1g_0^{-1},\ldots, g_0 t_r g_0^{-1},\\
&& g_0s\exp(\xi_0) g_0^{-1}=c)\in Q_1^{-1}(c_0).
\end{eqnarray*}

For $\ell>1$, define $\gamma_2:[0,1]\longrightarrow
G^{2\ell}\times \Ct \times G^2$ by
$$
\gamma_2(t)=(a_1(t), b_1(t),a_2(t),b_2(t), e,\ldots,e, d'_1(t),\ldots,d'_r(t), c_1(t),c_2(t)),
$$
where
\begin{eqnarray*}
a_1(t)&=&\tg_2(t) a \tg_2(t)^{-1}, \\
b_1(t)&=&\tg_2(t)\exp(-2t\xi_2')\tg_2(t)^{-1}, \\
a_2(t)&=&\tg_1(t) a \tg_1(t)^{-1}, \\
b_2(t)&=&\tg_1(t)\exp(-2t\xi_1')\tg_1(t)^{-1}, \\
d'_j(t)&=&\tg_1(t)t_j\tg_1(t)^{-1},\ \ j=1,\ldots,r,\\
c_1(t)&=& \tg_1(t)s\exp(t\xi_1)\tg_1(t)^{-1},\\
c_2(t)&=& \tg_2(t)\exp(t\xi_2)\tg_2(t)^{-1}.
\end{eqnarray*}
Then the image of $\gamma_2$ lies in $X_G^{\ell,r,2}(\Cr;z)$,
$$
\gamma_2(0)=(a,e,a,e,e,\ldots,e,t_1,\ldots,t_r,s,e)\in Q^{-1}(s,e),
$$
and
\begin{eqnarray*}
\gamma_2(1)&=&(g_2ag_2^{-1}, g_2\exp(-2\xi_2')g_2^{-1},
           g_1 a g_1^{-1}, g_1\exp(-2\xi_1') g_1^{-1},
       e,\ldots,e,\\
 & & g_1 t_1 g_1^{-1},\ldots,g_1t_rg_1^{-1}, c_1,c_2)\in Q^{-1}(c_1,c_2).\ \ \Box
\end{eqnarray*}

\bigskip

Recall that
$$
\Hom_{\Cr}(\pi_1(\lri),G)\cong X_G^{\ell,r,i}(\Cr;e)
$$
for $i=0,1,2$. So Proposition \ref{thm:connected} implies
that Theorem \ref{thm:irzero}--\ref{thm:nonzero}
hold when $G$ is simply connected:
\begin{cor}\label{thm:simply}
Let $G$ be a compact, connected, simply connected Lie group.
Let $\ell$ be a positive integer, let $r$ be a nonnegative integer,
let $\cC_1,\ldots,\cC_r$ be any conjugacy classes in $G$, and let
$z\in Z(G)$.  Then
$$
\Hom_{\Cr}(\pi_1(\slr),G)/G,\ \
\Hom_{\Cr}(\pi_1(\sone),G)/G
$$
are nonempty and connected;
$$
\Hom_{\Cr}(\pi_1(\stwo),G)/G,
$$
is nonempty, and is connected if $\ell>1$.
\end{cor}

\section{General case}\label{sec:general}

\subsection{Orientable surfaces}\label{sec:orientable}

Let $G, G_{ss}, H, D$ be as in Section \ref{sec:compact}, and
let $\pi:G\to G/G_{ss}=H/D$ be the projection. Then
$\pi$ descends to a map $\hat{\pi}: Con(G) \to H/D$. Let
$\pi_1:H\to H/D$ be the natural projection. We have the
following observation.
\begin{lm}\label{thm:image}
Let $G$ be a compact connected Lie group. Let $\ell$ be a positive integer.
Then
$$
\im(\Phi_{G,\ell,0})=G_{ss}.
$$
For $r>0$, we have
$$
\im(\Phi^{\Cr}_{G,\ell,r})= hG_{ss}
$$
where $h\in H$, $\pi_1(h)=\hat{\pi}(\cC_1)\cdots \hat{\pi}(\cC_r)\in H/D$.
\end{lm}

\paragraph{\bf Proof}
It is obvious from the definition that $\im\left(\Phi_{G,\ell,0}\right)\subset G_{ss}$.
Conversely, given $g\in G_{ss}$, let $\tg\in \tG_{ss}$ be a preimage of $g$
under $\rho_{ss}:\tG_{ss}\to G_{ss}$. By Fact \ref{thm:fiber}, there exist
$\ta,\tb\in \tG_{ss}$ such that $[\ta,\tb]=\tg$. Let $a=\rho_{ss}(\ta)$ and
$b=\rho_{ss}(\tb)$. Then
$$
g=[a,b]=\Phi_{G,\ell,0}(a,b,e,\ldots,e)\in \im\left(\Phi_{G,\ell,0}\right).
$$
So $G_{ss}\subset \im(\Phi_{G,\ell,0})$.

It follows from the above $r=0$ case that  $\im(\Phi^{\Cr}_{G,\ell,r})=G_{ss}$
if $\cC_j\subset G_{ss}$ for $j=1,\ldots, r$. In general, $\cC_j=G\cdot g_j$
for some $g_j\in G$, and $g_j=h_j g_j'$ for some $h_j\in H$ and $g_j'\in
G_{ss}$. We have $\cC_j= h_j\cC_j'$, where
$\cC_j'= G\cdot g_j'=G_{ss}\cdot g_j'\subset G_{ss}$, and
$\pi_1(h_j)=\hat{\pi}(\cC_j)$. So
\begin{equation}\label{eqn:CinGss}
\im(\Phi_{G,\ell,r}^{\cC_1',\ldots,\cC_r'})=G_{ss}.
\end{equation}
From the definition of $\Phi_{G,\ell,r}^{\cC_1,\ldots,\cC_r}$, one
sees that (\ref{eqn:CinGss}) implies
$$
\im(\Phi_{G,\ell,r}^{h_1\cC_1',\ldots,h_r\cC_r'})=h_1\cdots h_rG_{ss},
$$
or equivalently,
$$
\im(\Phi_{G,\ell,r}^{\Cr})=h G_{ss},
$$
where $h=h_1\cdots h_r\in H$ and $\pi_1(h_1\cdots
h_r)=\hat{\pi}(\cC_1)\cdots\hat{\pi}(\cC_r)$. $\ \ \Box$.

Let $\Cr$ be conjugacy classes of $G$. From the proof of
Lemma \ref{thm:image}, $\cC_j=h_j\cC_j'$ for some $h_j\in H$
and some conjugacy class $\cC_j'\subset G_{ss}$.
Recall that
$$
\Hom_{\Cr}(\pi_1(\Si^{\ell,r}_0),G)\cong X^{\ell,r,0}_G(\Cr;e)
=\left(\Phi^{\Cr}_{G,\ell,r}\right)^{-1}(e).
$$
By Lemma \ref{thm:image},
$\left(\Phi^{\Cr}_{G,\ell,r}\right)^{-1}(e)$ is nonempty iff $h_1\cdots
h_r\in D$. Note that if $d_1\cdots d_r\in G_{ss}$ for some
$(d_1,\ldots,d_r)\in \Ct$, then $h_1\cdots h_r\in D$ and
$d_1\cdots d_r\in G_{ss}$ for {\em all} $(d_1,\ldots,d_r)\in\Ct$.
So we have
\begin{lm}\label{thm:empty}
Let $G$ be a compact connected Lie group, and let
$\ell,r$ be positive integers. Then
$\Hom_{\Cr}(\pi_1(\Si^{\ell,r}_0,G))/G$ is nonempty iff
$d_1\cdots d_r\in G_{ss}$ for some (and therefore for all)
$(d_1,\dots,d_r)\in \Ct$.
\end{lm}

Let $\Cr$ be conjugacy classes in $G$ such that $d_1\cdots d_r\in G_{ss}$
for some $(d_1,\ldots,d_r)\in \Ct$.
From the above discussion, $\cC_j=h_j\cC_j'$ for some $h_j\in H$,
$\cC_j\subset G_{ss}$, and  $h=h_1\cdots h_r\in D$. Note that we
may replace $h_1$ and $\cC'_1$ by $h^{-1}h_1\in H$ and
$h\cC'_1\subset G_{ss}$, respectively, so we may assume that
$h_1\cdots h_r =e$. The diffeomorphism
$$
G^{\ell}\times\cC_1'\times\cdots\times\cC_r'
\longrightarrow G^{\ell}\times \Ct
$$
given by
$$
(\ab,d_1,\ldots,d_r)\mapsto (\ab, h_1d_1,\ldots,h_rd_r)
$$
induces an isomorphism of topological spaces
$$
X^{\ell,r,0}_G(\cC_1',\ldots,\cC'_r;e)\longrightarrow
X^{\ell,r,0}_G(\Cr;e).
$$
So
$$
\Hom_{\cC_1',\cdots,\cC_r'}(\pi_1(\Si^{\ell,r}_0),G)/G
\cong \Hom_{\Cr}(\pi_1(\Si^{\ell,r}_0),G)/G.
$$

From the above discussion, Theorem \ref{thm:izero} follows from Lemma \ref{thm:empty}
and Theorem \ref{thm:orientable}.
\begin{thm}\label{thm:orientable}
Let $G$ be a compact connected Lie group, and let $G_{ss}=[G,G]$ be the
maximal connected semisimple subgroup of $G$. Let $\ell, r$ be positive
integers, and let $\Cr$ be conjugacy classes of $G$ such that
$\cC_j\subset G_{ss}$ for $j=1,\ldots,r$. Then there is a bijection
$$
\pi_0(\Hom_{\Cr}(\pi_1(\Si^{\ell,r}_0),G))\longrightarrow \pi_1(G_{ss})/J_{\Cr}
$$
where $J_{\Cr}$ is defined as in Section \ref{sec:conjugate}.
\end{thm}
\paragraph{\bf Proof}
Let $\rho:\tG\to G$ and $\rho_{ss}:\tG_{ss}\to G_{ss}$ be  universal
coverings, as before. Let $\cD_j$ be a connected component of $\rho^{-1}(\cC_j)$.
Recall from Section \ref{sec:conjugate} that
$\cD_j$ is a conjugacy class of $\tG=\fh\times\tG_{ss}$, and
$\rho_j:\cD_j \to \cC_j$ is a finite cover of degree $|K_{\cD_j}|$,
where $K_{\cD_j}$ is a subgroup of
$\{0\}\times \Ker(\rho_{ss})\cong \pi_1(G_{ss})$.
$J_{\Cr}$ is the subgroup of $\{0\}\times \Ker(\rho_{ss})$ generated by $K_{\cD_1},\cdots, K_{\cD_r}$.
In the rest of this proof,  we will identify $\Ker(\rho_{ss})$ with $\{0\}\times \Ker(\rho_{ss})$.
Let
$$
P=\rho^{2\ell}\times \rho_1\times\cdots\times \rho_r:\tG^{2\ell}\times
\cD_1\times\cdots\times\cD_r
\to G^{2\ell}\times\Ct.
$$

Define
$$
o:\Hom_{\Cr}(\pi_1(\Si^{\ell,r}_0),G)\cong  X^{\ell,r,0}_G(\Cr;e)\to \Ker(\rho_{ss})/J_{\Cr}
$$
by
\begin{equation}\label{eqn:zero-o}
(a_1,b_1,\ldots,a_\ell,b_\ell,d_1,\ldots,d_r) \mapsto
\left[\prod_{i=1}^\ell [\ta_i,\tb_i]\prod_{j=1}^r \td_j \right]
\end{equation}
where
$$
(\ta_1,\tb_1,\ldots, \ta_\ell,\tb_\ell,
\td_1,\ldots,\td_r)\in P^{-1}(a_1,b_1,\ldots,a_\ell,b_\ell,d_1,\ldots,d_r),
$$
and $[k]$ denotes the image of $k\in \Ker(\rho_{ss})$ under the natural
projection $\Ker(\rho_{ss})\to \Ker(\rho_{ss})/J_{\Cr}$.
Note that the definition (\ref{eqn:zero-o}) does not depend on
the choice of $(\ta_1,\tb_1,\ldots,\ta_\ell, \tb_\ell,\td_1,\ldots,\td_r)$
because $\Ker\rho\subset Z(\tG)$ and
$$
\rho_j^{-1}(d_j)=\{ k\td_j\mid k\in K_{\cD_j}\}.
$$
The map $o$ descends to a continuous map
$$
\bar{o}:\Hom_{\Cr}(\pi_1(\Si^{\ell,r}_0),G)/G
\to\Ker(\rho_{ss})/J_{\Cr}.
$$

It is easy to check that $P$ restricts to a surjective map
$$
X_{\tG}^{\ell,r,0}(\cD_1,\ldots,\cD_r;k)
\cong \fh^{2\ell}\times X_{\tG_{ss}}^{\ell,r,0}(\cD_1,\ldots,\cD_r;k)
\longrightarrow o^{-1}([k])
$$
for any $k\in \Ker(\rho_{ss})$. By Fact \ref{thm:fiber},
$X_{\tG_{ss}}^{\ell,r,0}(\cD_1,\ldots,\cD_r;k)$
is nonempty and connected for any $k\in \Ker(\rho_{ss})$.
So $o^{-1}([k])$ and $\bar{o}^{-1}([k])$ are nonempty and connected
for each $[k]\in K/J_{\Cr}$. This completes the proof. $\ \ \Box$.

\bigskip

Finally, Theorem \ref{thm:irzero} follows immediately from Theorem
\ref{thm:orientable} because $\Hom(\pi_1(\Si^{\ell,0}_0),G)/G$ can
be identified with $\Hom_{\{e\}}(\pi_1(\Si^{\ell,1}_0),G)/G$ and
$J_{\{e\}}$ is trivial.

\subsection{Nonorientable surfaces}\label{sec:nonorientable}
Let $G, G_{ss}, H,D$ be defined as in Section \ref{sec:compact}, and
let $\Cr$ be $r>0$ conjugacy classes of $G$. We have seen in Section
\ref{sec:orientable} that $\cC_j=h_j\cC_j'$ for some $h_j\in H$ and
some conjugacy class $\cC_j'\subset G_{ss}$.
Let $s\in H$ be a square root of $h_1\cdots h_r$.
The diffeomorphism
\begin{eqnarray*}
 G^{2\ell}\times \Ct\times G &\longrightarrow
& G^{2\ell}\times \cC'_1\times\cdots\times\cC'_r\times G\\
(a_1,b_1,\ldots,a_\ell,b_\ell,d_1,\ldots,d_r,c) &\mapsto&
(a_1,b_1,\ldots,a_\ell,b_\ell,h_1^{-1} d_1,\ldots, h_r^{-1} d_r,sc)
\end{eqnarray*}
induces a homeomorphism of (topological) subspaces
$$
X_G^{\ell,r,1}(\Cr;e)\to X_G^{\ell,r,1}(\cC_1',\ldots,\cC_r';e).
$$
Similarly,
the diffeomorphism
\begin{eqnarray*}
 G^{2\ell}\times \Ct\times G^2
&\longrightarrow& G^{2\ell}\times \cC_1'\times\ldots\times\cC_r' \times G^2\\
(a_1,b_1,\ldots,a_\ell,b_\ell,d_1,\ldots,d_r,c_1,c_2) &\mapsto&
(a_1,b_1,\ldots,a_\ell,b_\ell,h_1^{-1} d_1,\ldots, h_r^{-1} d_r,c_1,sc_2)
\end{eqnarray*}
induces a homeomorphism
$$
X_G^{\ell,r,2}(\Cr;e)\to X_G^{\ell,r,2}(\cC'_1,\ldots,\cC'_r;e).
$$
So we have
$$
\Hom_{\Cr}(\pi_1(\Si^{\ell,r}_i),G)\cong
\Hom_{\cC_1',\cdots,\cC_r'}(\pi_1(\Si^{\ell,r}_i),G)
$$
for $\ell,r>0$ and $i=1,2$. Therefore, Theorem \ref{thm:nonzero}
follows from Theorem \ref{thm:nonorientable}.
\begin{thm}\label{thm:nonorientable}
Let $\rho:\tG\to G$ and $K=\Ker\rho\cong \pi_1(G)$ be defined as in Section \ref{sec:compact}.
Let $r$ be a positive integer, and let $\cC_1,\ldots,\cC_r$ be conjugacy
classes of $G$ such that $\cC_j\subset G_{ss}$ for $j=1,\ldots,r$.
Then for $\ell\geq i$, where $i=1,2$, there is a bijection
$$
\pi_0(\Hom_{\Cr}(\pi_1(\Si^{\ell,r}_i),G))\cong K/J'_{\Cr}
$$
where $J'_{\Cr}$ is the subgroup of $K$ generated by $J_{\Cr}$ and $2K=\{k^2\mid k\in K\}$.
\end{thm}

\begin{rem}
Notice that there is no condition on the conjugacy class for the moduli space to be nonempty.
\end{rem}

\paragraph{\bf Proof}
Let $\cD_j$ be a connected component of $\cC_j$. Since $\cC_j\subset G_{ss}$,
we may choose $\cD_j$ such that $\cD_j\subset \{0\}\times \tG_{ss}\subset \fh\times\tG_{ss}=\tG$.
Recall that $\rho_j:\cD_j\to \cC_j$ is a finite cover with degree $|K_{\cD_j}|$, and
$J_{\Cr}$ is the subgroup of $K$ generated by $K_{\cD_1},\cdots,K_{\cD_r}$.
Set
\begin{eqnarray*}
&P_1=\rho^{2\ell}\times \rho_1\times\cdots\rho_r\times \rho:&\tG^{2\ell}\times
\cD_1\times\cdots\times\cD_r\times \tG
\to G^{2\ell}\times\Ct\times G\\
&P_2=\rho^{2\ell}\times \rho_1\times\cdots\rho_r\times \rho^2:&\tG^{2\ell}\times
\cD_1\times\cdots\times\cD_r\times \tG^2
\to G^{2\ell}\times\Ct\times G^2
\end{eqnarray*}
For $i=1,2$ we define
$$
o:X^{\ell,r,i}_G(\Cr;e)\to K/J'_{\Cr}
$$
as follows. For $i=1$, $o$ is given by
\begin{equation}\label{eqn:one-o}
(a_1,b_1,\ldots,a_\ell,b_\ell,d_1,\ldots,d_r,c) \mapsto
\left[\prod_{i=1}^\ell [\ta_i,\tb_i]\prod_{j=1}^r \td_j \tc^2\right]
\end{equation}
where
$$
(\ta_1\tb_1,\ldots,\ta_\ell,\tb_\ell,
\td_1,\ldots,\td_r,\tc) \in P_1^{-1}(a_1,b_1,\ldots,a_\ell,b_\ell,d_1,\ldots,d_r,c),
$$
and $[k]$ denotes the image of $k$ under the natural projection $K\to
K/J'_{\Cr}$. For $i=2$, $o$ is given by
\begin{equation}\label{eqn:two-o}
(a_1,b_1,\ldots,a_\ell,b_\ell,d_1,\ldots,d_r,c_1,c_2) \mapsto
\left[\prod_{i=1}^\ell [\ta_i,\tb_i]\prod_{j=1}^r \td_j\tc_1^2\tc_2^2\right]
\end{equation}
where
$$
(\ta_1\tb_1,\ldots,\ta_\ell,\tb_\ell, \td_1,\ldots,\td_r,\tc_1,\tc_2) \in
P_2^{-1}(a_1,b_1,\ldots,a_\ell,b_\ell,d_1,\ldots,d_r,c_1,c_2),
$$
Note that the definitions (\ref{eqn:one-o}), (\ref{eqn:two-o}) do not
depend on the choices of
$$
(\ta_1,\tb_1,\ldots,\ta_\ell,\tb_\ell,\td_1,\ldots,\td_r,\tc),\ \
(\ta_1,\tb_1,\ldots,\ta_\ell, \tb_\ell,\td_1,\ldots,\td_r,\tc_1,\tc_2)
$$
respectively. The map $o$ descends to a continuous map
$$
\bar{o}:\Hom_{\Cr}(\pi_1(\Si^{\ell,r}_i,G)/G\to K/J'_{\Cr}
$$
where $i=1,2$.

For $i=1,2$, it is easy to check that $P_i$ restricts to  a surjective map
$$
X_{\tG}^{\ell,r,i}(\cD_1,\ldots,\cD_r;(X,z))
\cong \fh^{2\ell+i-1}\times X_{\tG_{ss}}^{\ell,r,i}(\cD_1,\ldots,\cD_r;z)
\longrightarrow o^{-1}([(X,z)])
$$
for any $(X,z)\in K\subset H\times Z(\tG_{ss})$. By Fact \ref{thm:fiber},
$X_{\tG_{ss}}^{\ell,r,i}(\cD_1,\ldots,\cD_r;z)$
is nonempty and connected if $\ell\geq i$. So $o^{-1}([k])$ is
nonempty and connected if $\ell\geq i$. This completes the proof.  $\ \ \Box$

\bigskip

Finally, Theorem \ref{thm:rzero} follows directly from Theorem
\ref{thm:nonorientable} because $\Hom(\pi_1(\Si^{\ell,0}_i),G)/G$
can be identified with $\Hom_{\{e\}}(\pi_1(\Si^{\ell,1}_i),G)/G$
and $J_{\{e\}}$ is trivial.

\section{Moduli spaces of flat principal $G$-bundles} \label{sec:moduli}

Let $G$ be a compact connected Lie group, and let $\Si$ be a compact
surface. Recall from Section \ref{sec:variety}
that any compact surface is of the form $\Si^{\ell,r}_{i}$
for some nonnegative integers $\ell,r$ and $i=0,1,2$.

A principal $G$-bundle $P$ over a compact surface $\Si$ is trivial
iff it admits a cross section. The obstruction classes $O_n(P)\in H^n(\Si;\pi_{n-1}(G))$ are
obstructions to the existence of a cross section of $P$. They are topological
invariants of $P$. For a surface, the only obstructions are $O_1$ and
$O_2$. Here we consider connected Lie groups so the first obstruction class
$O_1$ is trivial and we are left with the second obstruction class $O_2$
only.

Let $\Prin_G(\Si)$ denote the moduli space of topological principal $G$
bundles over $\Si$. Then $P\mapsto O_2(P)$ defines a bijection
$$
O_2:\Prin_G(\Si)\to H^2(\Si;\pi_1(G))
$$
where
$$
H_2(\lri;\pi_1(G))=\left\{\begin{array}{ll}
\pi_1(G), &  r=i=0,\\
\pi_1(G)/2\pi_1(G), & r=0, i=1,2,\\
0 & r>0.
\end{array}\right.
$$

We first consider the case  $r>0$. From the above discussion,
every principal $G$-bundle over $\lri$ is trivial.
The space of connections on the trivial bundle $\lri\times G$ can be
identified with $\Omega^1(\lri,\fg)$, the space of $\fg$-valued 1-forms on $\lri$.
Let $\Hol_j:\Omega^1(\lri,\fg)\to G$ be the holonomy around the boundary $B_j$.
Given conjugacy classes $\Cr$ of $G$, let
$$
\mlri=\frac{\{A\in\Omega^1(\lri,\fg)\mid F_A=0, \Hol_j(A)\in\cC_j\textup{ for }
  j=1,\ldots, r\} }{C^\infty(\lri,G)}
$$
be the moduli space of gauge equivalence classes of flat connections
on $\lri\times G$ with holonomy around $B_j$ in $\cC_j$.
With suitable choices of orientation on $B_j$, we have
$$
\mlri\cong \Hom_{\Cr}(\pi_1(\Si^{\ell,r}_i),G)/G.
$$
So in Theorem \ref{thm:izero} and Theorem \ref{thm:nonzero},
we may replace $\holC/G$ by $\cM_G(\Si;\Cr)$.

We next consider the case $r=0$.  A principal $G$-bundle on $\Si^{\ell,0}_i$
may be topologically nontrivial. Let $\cM_G(\Si^{\ell,0}_i)$ be the moduli
space of gauge equivalence classes of flat $G$-bundles on $\Si^{\ell,0}_i$.
Here a flat $G$-bundle is a principal $G$-bundle together with a flat $G$
connection.
Then
$$
\cM_G(\Si^{\ell,0}_i)\cong \Hom(\pi_1(\Si^{\ell,0}_i), G)/G.
$$
Let
$$
F: \Hom(\pi_1(\Si^{\ell,0}_i),G)/G\to \Prin_G(\Si^{\ell,0}_i),
$$
be the map which sends a flat $G$-bundle to its underlying topological
principal $G$-bundle. The discussion in \cite[Appendix A]{HL1} shows that
the map
$$
O_2\circ F:\Hom(\pi_1(\Si^{\ell,0}_i),G)/G\to H^2(\Si^{\ell,0}_i;\pi_1(G))
$$
coincides with
$$
\bar{o}:\Hom_{\{e\}}(\pi_1(\Si^{\ell,1}_i),G)/G\cong\Hom(\pi_1(\Si^{\ell,0}_i),G)/G
\to H^2(\Si^{\ell,0}_i;\pi_1(G))$$ defined in the proofs of Theorem \ref{thm:orientable}
and Theorem \ref{thm:nonorientable} in Section \ref{sec:general}. See \cite[Appendix A]{HL1}
for various interpretations of the obstruction map.

Given a topological principal $G$-bundle $P$ over $\Si$, let
$\cM(P)$ denote the moduli space of flat $G$-connections on $P$.
Note that the subgroup $\pi_1(G_{ss})$ of $\pi_1(G)$
consists of torsion elements in $\pi_1(G)$. From the above discussion, our proof
of Theorem \ref{thm:orientable} gives the following statement:
\begin{tm}\label{thm:orientableM}
Let $\Si$ be a connected, closed, compact, orientable surface with
genus $\ell>0$. Let $G$ be a compact connected Lie group, and let
$P$ be a principal $G$-bundle over $\Si$. Then $\cM(P)$ is
nonempty if and only if the obstruction class $O_2(P)\in
H^2(\Si;\pi_1(G))\cong \pi_1(G)$ is a torsion element. In this
case, $\cM(P)$ is connected.
\end{tm}

Similarly, our proof of Theorem \ref{thm:nonorientable} gives the following statement:
\begin{tm}\label{thm:nonorientableM}
Let $\Si$ be a closed, compact, orientable surface which
is homeomorphic to $m$ copies of the real projective plane.
Let $G$ be a compact connected Lie group, and let $P$ be a
principal $G$-bundle over $\Si$. Then $\cM(P)$ is nonempty.
Moreover, $\cM(P)$ is connected if $m\neq 1,2,4$.
\end{tm}

Note that Theorem \ref{thm:orientableM} implies Theorem \ref{thm:irzero},
and Theorem \ref{thm:nonorientableM} implies Theorem \ref{thm:rzero}.
We now outline a proof of Theorem \ref{thm:orientableM} extracted from
\cite{ym}. Let $\Si$, $G$ be as in Theorem \ref{thm:orientableM}, and let
$P$ be a principal $G$-bundle over $\Si$. By \cite[Proposition 6.16]{ym},
$P$ possesses a central Yang-Mills connection. It is proved in
\cite[Section 10]{ym} that the semi-stable stratum, which
is an open and nonempty subset in the space of all $G$-connections on $P$,
contains a unique connected component of the space of Yang-Mills connections
on $P$. This connected component consists of central Yang-Mills
connections on $P$. It is proved in \cite[Section 12]{ym} that central Yang-Mills
connections on $P$ achieve the absolute minimum of Yang-Mills functional on
$P$. This absolute minimum is a topological invariant $L(P)$ of $P$ which
vanishes iff $O_2(P)\in\pi_1(G)$ is a torsion element.
The $G$-bundle $P$ admits a flat $G$-connection iff $L(P)=0$,
and in this case, the space of flat $G$-connections on $P$ is the space
of central Yang-Mills connections on $P$.  This proves Theorem
\ref{thm:orientableM}.

\paragraph{\bf Acknowledgments:}
It is a pleasure to thank L. Jeffrey, E. Meinrenken, and E. Xia for their suggestions and help
during the preparation of this note.


\begin{thebibliography}{AMW}

\bibitem[AB]{ym} M.F. Atiyah and R. Bott,
{\em Yang-Mills equations over Riemann surfaces}, Philos. Trans.
Roy. Soc. London Ser. A \textbf{308} (1983), no. 1505, 523--615.

%\bibitem[ABS]{spin} M.F. Atiyah, R. Bott, and A. Shapiro,
%{\em Clifford modules}, Topology \textbf{3} (1964) suppl. 1,
%3--38.

\bibitem[AMM]{gvm} A. Alekseev, A. Malkin, and E. Meinrenken,
 {\em Lie group valued moment maps},
 J. Differential Geom. \textbf{48} (1998), no. 3, 445--495.

\bibitem[AMW]{dhm} A. Alekseev, E. Meinrenken, and C. Woodward,
 {\em Duistermaat-Heckman measures and moduli spaces of flat bundles over
 surfaces},
 Geom. Funct. Anal. \textbf{12} (2002), no. 1, 1--31.

\bibitem[BD]{bd} Theodor Br\"{o}cker and Tammo tom Dieck,
 {\em Representations of compact Lie groups},
 Graduate Texts in Mathematics, 98. Springer-Verlag, New York, 1985.

\bibitem[HM]{HM} K.H. Hofmann and S.A. Morris,
{\em The structure of compact groups.
A primer for the student---a handbook for the expert},
de Gruyter Studies in Mathematics, 25.
Walter de Gruyter \& Co., Berlin, 1998.

\bibitem[G1]{G1} W.M. Goldman,
 {\em The symplectic nature of fundamental groups of surfaces},
 Adv. in Math. \textbf{54} (1984), no. 2, 200--225.

\bibitem[G2]{G2} W.M. Goldman,
{\em Topological components of spaces of representations}, Invent.
Math. \textbf{93} (1988), no. 3, 557--607.

\bibitem[HL1]{HL1} N.-K. Ho and C.-C.M. Liu,
{\em Connected Components of the Space of Surface Group Representations}
IMRN (2003), no. 44, 2359--2372.

\bibitem[HL2]{HL2} N.-K. Ho and C.-C.M. Liu,
{\em On the Connectedness of Moduli Spaces of Flat
Connections over Compact Surfaces}, arXiv:math.SG/0211388,
to appear in Canadian Journal of Mathematics.

\bibitem[Hu]{Hu} J.E. Humphreys,
{\em Reflection Groups and Coxeter Groups}, Cambridge Studies in Advanced
Mathematics, vol.29, Cambridge University Press, Cambridge, 1990.

\bibitem[K]{K} A.W. Knapp,
{\em Lie Groups Beyond an Introduction},
Progress in Mathematics, 140, Birkh\"{a}user Boston, Inc., Boston, MA, 1996.

\bibitem[Li]{li} J. Li,
 {\em The space of surface group representations},
 Manuscripta Math. \textbf{78} (1993), no. 3, 223--243.

\bibitem[Ra]{Ra} A. Ramanathan,
{\em Stable principal bundles on a compact Riemann surface},
Math. Ann.  \textbf{213}  (1975), 129--152.

\end{thebibliography}
\end{document}